\RequirePackage[l2tabu, orthodox]{nag}
\RequirePackage{snapshot}

\documentclass[10pt,twocolumn]{article}

\makeatletter
\if@twocolumn
  \usepackage[dvips,letterpaper,top=0.5in, bottom=0.5in, left=0.75in, right=0.5in,includefoot,heightrounded]{geometry}
\else
  \usepackage[dvips,letterpaper,margin=1in,includefoot,heightrounded]{geometry}
\fi

\usepackage{srcltx}

\usepackage{amsmath}
\usepackage{amssymb,amsfonts}

\usepackage{abstract}

\usepackage[usenames,dvipsnames]{color}
\usepackage[usenames,dvipsnames]{xcolor}
\usepackage{subfigure}

\usepackage{booktabs}

\usepackage{setspace}
\usepackage{flushend}
\usepackage{multicol}

\usepackage{cite}
\usepackage{url}
\usepackage[normalem]{ulem}

\usepackage{enumerate}

\usepackage{hyperref}

\usepackage[sc,tiny]{titlesec}
       \usepackage{mathptmx}
       \usepackage[mathscr]{euscript}

\newtheorem{proposition}{Proposition}

\newtheorem{corollary}{Corollary}

\def\QED{\mbox{$\square$}}
\def\proof{\noindent\textit{Proof:~}}
\def\endproof{\hspace*{\fill}~\QED\par\endtrivlist\unskip}



\title{%
Fourier Eigenfunctions, Uncertainty Gabor Principle and
Isoresolution~Wavelets
}

\author{%
L. R. Soares%
\thanks{%
L. R. Soares
was with the
Grupo de Pesquisa em Comunica\c{c}\~oes (\textsc{codec}),
Departamento de Eletr\^onica e Sistemas,
Universidade Federal de Pernambuco.
}
\and
H. M. de Oliveira%
\thanks{%
H. M. de Oliveira
was with the
Grupo de Pesquisa em Comunica\c{c}\~oes (\textsc{codec}),
Departamento de Eletr\^onica e Sistemas,
Universidade Federal de Pernambuco.
Currently he is with the
Signal Processing Group,
Departamento de Estat\'{\i}stica,
Universidade Federal de Pernambuco.
Email:~\url{hmo@ufpe.br}
}
\and
R.~J.~Cintra%
\thanks{%
R.~J.~Cintra
was
with the Graduate Program in Electrical Engineering,
Universidade Federal de Pernambuco.
Currently he is
with
the Signal Processing Group,
Departamento de Estat\'{\i}stica,
Universidade Federal de Pernambuco.
Email:~\url{rjdsc@de.ufpe.br}
}
\and
R. M. Campello de Souza%
\thanks{%
R. M. Campello de Souza
is with
Departamento de Eletr\^onica e Sistemas,
Universidade Federal de Pernambuco.
Email:~\url{ricardo@ufpe.br}
}
}

\date{}

\newcommand{\myabstract}{%
Shape-invariant signals under Fourier transform are investigated leading to a class of eigenfunctions for the Fourier operator. The classical uncertainty Gabor-Heisenberg principle is revisited and the concept of isoresolution in joint time-frequency analysis is introduced. It is shown that any Fourier eigenfunction achieve isoresolution. It is shown that an isoresolution wavelet can be derived from each known wavelet family by a suitable scaling.
}

\newcommand{\mykeywords}{%
Gabor-Heisenberg inequality, Fourier eigenfunctions, Isoresolution wavelets, time-frequency analysis.
}

\begin{document}

\makeatletter
\if@twocolumn

\twocolumn[%
  \maketitle
  \begin{onecolabstract}
    \myabstract
  \end{onecolabstract}
  \begin{center}
    \small
    \textbf{Keywords}
    \\\medskip
    \mykeywords
  \end{center}
  \bigskip
]
\saythanks

\else

  \maketitle
  \begin{abstract}
    \myabstract
  \end{abstract}
  \begin{center}
    \small
    \textbf{Keywords}
    \\\medskip
    \mykeywords
  \end{center}
  \bigskip
  \onehalfspacing
\fi

\section{Preliminaries}

The Fourier transform is often interpreted as a linear operator~$\mathcal{F}$.
An interesting problem in this framework is to find out
the eigenfunctions
in the language of operators~\cite{HERS64,SOKORED66,PEIDIN02}.
Let $\mathcal{V}$ be a vector space equipped with a linear transform,
$T: \mathcal{V}\to\mathcal{V}$, $\mathbf{v} \mapsto T(\mathbf{v})$.
Under the linear transform~$T$,
eigenfunctions are solutions of
$T(\mathbf{v})= \lambda \cdot \mathbf{v}$,
which corresponds here to
$\operatorname{\mathcal{F}}\{ f(t) \}(\omega) = \lambda \cdot f(\omega)$
where
$f \in L^2(\mathbb{R})$
and
$\lambda$ is a scalar.
They are a quite remarkable class of functions, which preserves the shape under Fourier transform: Both the signal and its spectrum (time and frequency representation) have the same shape.
In joint time-frequency representation~\cite{COH95,QIACHE99} this feature can represent a very good balance between the two domains. It is well known that the Gaussian pulse is a signal whose shape is preserved under the Fourier operator:
\begin{align*}
e^{-t^2/2}
\stackrel{\mathcal{F}}{\longleftrightarrow}
\sqrt{2\pi}
\cdot
e^{-\omega^2/2}
.
\end{align*}

This can easily be derived by writing
\begin{align*}
\frac{1}{\sqrt{2\pi}}
\int_{-\infty}^\infty
e^{-t^2/2}
\cdot
e^{j \omega t}
\operatorname{d}t
=
F(\omega)
.
\end{align*}
Deriving this equation and using integral by parts,
one notice that:
$\frac{\operatorname{d}}{\operatorname{d}\omega}F(\omega)=-\omega F(\omega)$.
The solution of the differential equation
$\frac{\operatorname{d}}{\operatorname{d}\omega}F(\omega)
+\omega F(\omega)
=
0$
under the initial condition $F(0)=1$ is
$F(\omega) = e^{-\omega^2/2}$.
It follows promptly that
$\lambda = \sqrt{2\pi}$.

The question is: Are there other eigenfunctions?
This matter is addressed in the next section.
It is worthwhile to bear in mind that some results in this paper are deliberately \emph{non nova, sed nove}.

\section{Shape-invariant Signals: Eigenfunctions of the Fourier Operator}

Let $\mathscr{E}\{\cdot\}$ and
$\mathscr{O}\{\cdot\}$
denote the functionals that extract the even and
odd part of a given signal,
respectively.

\begin{proposition}
\label{proposition-1}
Let
$f(t) \stackrel{\mathcal{F}}{\longleftrightarrow}  F(\omega)$
be an arbitrary Fourier transform pair.
Then the signal
\begin{align*}
h(t)
=
\sqrt{2\pi}
\cdot
\mathscr{E}
\{ f(t) \}
+
\mathscr{E}
\{ F(\omega) \}
\end{align*}
is invariant under the Fourier transform.
Furthermore,
we have that:
$
H(\omega)
=
\operatorname{\mathcal{F}}
\{ h(t) \}(\omega)
=
\sqrt{2\pi}
\cdot
h(\omega)$.
\end{proposition}
\proof
It follows from the definition of
$h(\cdot)$
that
\begin{align*}
2 \cdot h(t)
=
\sqrt{2\pi}
\cdot
\left[
f(t) + f(-t)
\right]
+
\left[
F(t) + F(-t)
\right]
.
\end{align*}
Taking the Fourier transform,
\begin{align*}
2 \cdot H(\omega)
=
\sqrt{2\pi}
\cdot
\left[
F(\omega) + F(-\omega)
\right]
+
\left[
2\pi f(-\omega) + 2\pi f(\omega)
\right]
.
\end{align*}
and the proof follows.
\endproof

\begin{corollary}
Each even function
$f(t) \stackrel{\mathcal{F}}{\longleftrightarrow}  F(\omega)$
induces a Fourier invariant
$h(t) = \sqrt{2\pi} f(t) + F(t)$.
\end{corollary}

For instance, the following signals
\begin{align*}
h_1(t)
&=
\sqrt{2\pi}
\cdot
\frac{1}{1+t^2}
+
\pi
e^{-|t|}
,
\\
h_2(t)
&=
\sqrt{2\pi} |t|
-
\frac{2}{t^2}
\end{align*}
have spectra with similar shape.
Another remarkable example is:
\begin{align}
\label{equation-4}
\operatorname{sech}
\left(
\sqrt{\frac{\pi}{2}}
t
\right)
\stackrel{\mathcal{F}}{\longleftrightarrow}
\sqrt{2\pi}
\operatorname{sech}
\left(
\sqrt{\frac{\pi}{2}}
\omega
\right)
,
\end{align}
where
$\operatorname{sech}(\cdot)$
is the hyperbolic secant function.

\begin{proposition}
Let
$f(t) \stackrel{\mathcal{F}}{\longleftrightarrow}  F(\omega)$
be an arbitrary Fourier transform pair.
Then the signal
\begin{align*}
h(t)
=
\sqrt{2\pi}
\cdot
\mathscr{O}
\left\{
f(t)
\right\}
-
\mathscr{O}
\left\{
F(t)
\right\}
\end{align*}
is an invariant under Fourier transform.
Furthermore,
$\operatorname{\mathcal{F}}\left\{ h(t) \right\} = -\sqrt{2\pi} h(\omega)$.
\end{proposition}
\proof
The proof is similar to the proof of Proposition~\ref{proposition-1}.

\begin{corollary}
Each odd function
$f(t) \stackrel{\mathcal{F}}{\longleftrightarrow}  F(\omega)$
induces a Fourier invariant
$
h(t)
=
\sqrt{2\pi}
f(t)
-
F(t)
$.
\end{corollary}

Let us now focus on a particular and important class of Fourier invariant, which generates an orthogonal and complete set.
To begin with, let us denote by
$\mathcal{E}$
the class of eigenfunctions of the Fourier operator defined according to the following proposition.

\begin{proposition}
A signal $f(t)$ is in $\mathcal{E}$
if, and only if,
the signal
$f$ satisfies the differential equation
$
\frac{\operatorname{d}^2}{\operatorname{d}t^2}
f(t)
-
t^2
\cdot
f(t)
=
\kappa
\cdot
f(t)$,
for some scalar
$\kappa \in \mathbb{C}$.
\end{proposition}
\proof
We begin demonstrating the suffiency.
By hypothesis,
we have that
\begin{align*}
f(t)
\stackrel{\mathcal{F}}{\longleftrightarrow}
\lambda
f(\omega)
.
\end{align*}

The properties of time and frequency differentiation for 
$\mathcal{F}$
give:
\begin{align*}
\frac{\operatorname{d}^2}{\operatorname{d}t^2}
f(t)
&
\stackrel{\mathcal{F}}{\longleftrightarrow}
(j\omega)^2
\lambda
f(\omega),
\\
(-jt)^2
f(t)
&
\stackrel{\mathcal{F}}{\longleftrightarrow}
\lambda
\frac{\operatorname{d}^2}{\operatorname{d}\omega^2}
f(\omega)
.
\end{align*}
Adding above expressions%
\footnote{N.B. Subtracting: 
$\frac{\operatorname{d}^2}{\operatorname{d}t^2}
f(t)
+ t^2 f(t)
\stackrel{\mathcal{F}}{\longleftrightarrow}
-
\lambda
\left[
\frac{\operatorname{d}^2}{\operatorname{d}\omega^2}
f(\omega)
+
\omega^2
f(\omega)
\right]$.},
we derive
\begin{align*}
\frac{\operatorname{d}^2}{\operatorname{d}t^2}
f(t)
- t^2 f(t)
\stackrel{\mathcal{F}}{\longleftrightarrow}
-
\lambda
\left[
\frac{\operatorname{d}^2}{\operatorname{d}\omega^2}
f(\omega)
-
\omega^2
f(\omega)
\right]
.
\end{align*}

Thus,
the signal
$\frac{\operatorname{d}^2}{\operatorname{d}t^2}f(t)- t^2 f(t)$
has also its shape preserved, provided that $f$
itself preserves its shape.
Therefore, 
$\frac{\operatorname{d}^2}{\operatorname{d}t^2}f(t)- t^2 f(t)\in\mathcal{E}$, 
that is, we are looking for signals such that
$\frac{\operatorname{d}^2}{\operatorname{d}t^2}f(t)- t^2 f(t)=\kappa f(t)$,
since they have identical eigenvalues.

Now we demonstrate the necessity.

By hypohthesis,
the signal $f(t)$
satisfies the differential equation 
$\frac{\operatorname{d}^2}{\operatorname{d}t^2}f(t)- t^2 f(t)=\kappa f(t)$,
$\kappa \in \mathbb{C}$.
Applying the operator~$\mathcal{F}$,
we obtain: 
\begin{align*}
(j\omega)^2
F(\omega)
+
\frac{\operatorname{d}^2}{\operatorname{d}\omega^2}
F(\omega)
=
\kappa
\lambda
F(\omega)
.
\end{align*}
Thus,
$
\frac{\operatorname{d}^2}{\operatorname{d}\omega^2}
F(\omega)
-\omega^2
F(\omega)=
\kappa \lambda F(\omega)$, i.e., 
its spectrum also obeys a similar differential equation.
Therefore, $f$ and $F$ have identical shape,
since they are solutions of the same differential equation.
\endproof

The key equation for shape-invariant signal is
thus
$\frac{\operatorname{d}^2}{\operatorname{d}t^2}f(t)- t^2 f(t)=\kappa f(t)$.
Let us try solutions of the form 
\begin{align*}
f(t)
=
p(t)
e^{-t^2/2}
,
\end{align*}
where $p(t)$ is a function to be determined.
Therefore, 
\begin{align*}
\frac{\operatorname{d}^2}{\operatorname{d}t^2}
\left[
p(t)
e^{-t^2/2}
\right]
-
t^2
p(t)
e^{-t^2/2}
=
\kappa
p(t)
e^{-t^2/2}
.
\end{align*}
After simple algebraic manipulations,
we derive
\begin{align*}
\frac{\operatorname{d}^2}{\operatorname{d}t^2}
p(t)
-
2t
\frac{\operatorname{d}}{\operatorname{d}t}
p(t)
+
(\kappa + 1)
p(t)
=
0
,
\end{align*}
where $n$ is a integer.

A standard differential equation of the above form~\cite{ABRASTE68}
is
\begin{align}
\label{equation-11}
\frac{\operatorname{d}^2}{\operatorname{d}t^2}
p(t)
-
2t
\frac{\operatorname{d}}{\operatorname{d}t}
p(t)
+
2n
p(t)
=
0
,
\end{align}
where $n$ is a integer.
Thus, for a suitable choice 
$\kappa = -(2n+1)$
(eigenvalues),
the solutions $p(t)$
are exactly Hermite polynomials~\cite{ABRASTE68}, 
which form a complete orthogonal system.
Thus, we have:
\begin{align*}
p(t)
=
H_n(t),
\end{align*}
where
\begin{align*}
H_0(t) &=1,
\\
H_1(t) &=2t,
\\
H_2(t) &=-2+4t^2,
\\
H_3(t) &=-12t+8t^3,
\\
H_4(t) &= 12-48t^2+16t^4,
\\
&\vdots
\end{align*}

\begin{proposition}
Possible eigenvalues of the Fourier transform
are the four roots of the unit ($\pm1, \pm j$) 
times
$\sqrt{2\pi}$.
\end{proposition}
\proof
Let us denote by
$\mathcal{F}^{(n)}$
the operator corresponding to iterate $n$ times 
the operator~$\mathcal{F}$.
Let
$t
\stackrel{\mathcal{F}}{\longleftrightarrow}
\omega
\stackrel{\mathcal{F}}{\longleftrightarrow}
\omega'
\stackrel{\mathcal{F}}{\longleftrightarrow}
\Omega$
be the Fourier domain variables for the iterate Fourier transform.
Observe that,
for $f\in \mathcal{E}$,
we have:
\begin{equation}
\label{equation-13}
\begin{split}
\operatorname{\mathcal{F}}^{(2)}
\left\{
f(t)
\right\}
(\omega')
&=
2\pi
f(-\omega')
,
\\
\operatorname{\mathcal{F}}^{(4)}
\left\{
f(t)
\right\}
(\Omega)
&=
4\pi^2
f(\Omega)
.
\end{split}
\end{equation}
But,
\begin{equation}
\label{equation-14}
\begin{split}
\operatorname{\mathcal{F}}^{(2)}
\left\{
f(t)
\right\}
(\omega')
&=
\lambda^2
f(-\omega')
,
\\
\operatorname{\mathcal{F}}^{(4)}
\left\{
f(t)
\right\}
(\Omega)
&=
\lambda^4
f(\Omega)
.
\end{split}
\end{equation}
From~\eqref{equation-13} and~\eqref{equation-14},
it follows that $\lambda / \sqrt{2\pi} \in \mathbf{C}$
has order~4.
\endproof

We conclude that 
$\left\{
\psi_n(t)
=
H_n(t)
e^{-t^2/2}
\right\}_{n=0}^\infty$ 
are shape-invariant under Fourier operator associated to
$\lambda_n = (-j)^n \sqrt{2\pi}$. 
Therefore,
\begin{align}
\label{equation-15}
H_n(t)
e^{-t^2/2}
\stackrel{\mathcal{F}}{\longleftrightarrow}
(-j)^n
\sqrt{2\pi}
H_n(\omega)
e^{-\omega^2/2}
.
\end{align}
Another interpretation can be derived evoking
Rodrigues' formula~\cite{ABRASTE68}:
\begin{align*}
H_n(t)
=
(-1)^n
e^{t^2}
\frac{\operatorname{d}^n}{\operatorname{d}t^n}
e^{-t^2}
.
\end{align*}
The 2nd-order differential equation hold by invariant signals is 
\begin{align*}
\frac{\operatorname{d}^2}{\operatorname{d}x^2}
y
+
(2n+1-x^2)
y
=
0
.
\end{align*}
The above differential equation is exactly the celebrated Schr\"odinger equation for the harmonic oscillator~\cite{BEI67}.

\section{Consequences on the Time-Frequency Plane}

Let us now investigate certain consequences of eigenfunctions of the Fourier operator on the 
time-frequency plane~\cite{COH95,OLIBAR00}. 

Let $f(t)$ be a finite energy signal $E$, 
equipped with Fourier transform, $F(\omega)$.
The time and frequency moments of $f$ are defined by:
\begin{align*}
\overline{t^n}
&
=
\frac
{
\int_{-\infty}^\infty
f^\ast(t)
t^n
f(t)
\operatorname{d}t
}
{
\int_{-\infty}^\infty
f^\ast(t)
f(t)
\operatorname{d}t
}
,
\\
&
=
\frac{1}{E}
\int_{-\infty}^\infty
t^n
|f(t)|^2
\operatorname{d}t
\\
\overline{\omega^n}
&
=
\frac
{
\int_{-\infty}^\infty
F^\ast(\omega)
\omega^n
F(\omega)
\operatorname{d}\omega
}
{
\int_{-\infty}^\infty
F^\ast(\omega)
F(\omega)
\operatorname{d}\omega
}
\\
&
=
\frac{1}{2\pi E}
\int_{-\infty}^\infty
\omega^n
|F(\omega)|^2
\operatorname{d}\omega
.
\end{align*}

By analogy to Probability Theory, 
the term $|f(t)|^2/E$ denotes a ``time-domain'' energy density, 
where $E$ is a normalising factor so as to make the whole integral of the density be equal to the unity.
It is customary to deal with the energy spectral density 
$G(\omega)=|F(\omega)|^2$,
whose integral over a frequency band gives the energy content of the signal within such a band.
Let us suppose in the sequel, without loss of generality, that $E=1$ (energy normalised signals).

The ``effective duration''
(respectively the ``effective frequency width'') 
of a signal $f(t)$ (respectively $F(\omega)$) 
is defined according to:
\begin{align*}
\Delta t
&=
\sqrt{2 \pi \overline{(t-\overline{t})^2}}
\quad
\text{r.m.s duration,}
\\
\Delta f
&=
\sqrt{2 \pi \overline{(f-\overline{f})^2}}
\quad
\text{r.m.s bandwidth}
,
\end{align*}
where
$\Delta t$ 
and
$\Delta f$
correspond to the standard deviation, i.e., spreading measures.
However, other common and much handier definitions are 
\begin{align*}
\Delta_t 
&
=
\sqrt{\overline{(t-\overline{t})^2}}
,
\\
\Delta_\omega
&
=
\sqrt{\overline{(f-\overline{f})^2}}
.
\end{align*}
Clearly,
$\Delta_t = \Delta t / \sqrt{2\pi}$
and
$\Delta_\omega = \sqrt{2\pi} \Delta f$.

\subsection{Revisiting the Gabor Principle}

By applying arguments from quantum mechanics~\cite{BEI67},
Gabor~\cite{GAB46, GAB53} 
derived an uncertainty relation nowadays called Gabor-Heisenberg principle for signals:
$\Delta t \cdot \Delta f \geq 1/2$,
proving that time and frequency cannot be exactly measured (simultaneously). 
The Gabor-Heisenberg uncertainty principle states a lower bound on the product $\Delta t \cdot \Delta\omega$,
or alternatively:
\begin{align}
\label{equation-20}
\Delta_t\cdot\Delta_\omega \geq 1/2
.
\end{align}

\begin{proposition}
The Gabor lower bound is only achieved by the first invariant signal (eigenfunctions of $\mathcal{F}$ operator).
\end{proposition}
\noindent
\emph{Sketch of the proof:}
From~\eqref{equation-20},
the bound is achieved
if, and only if,
$\frac{\operatorname{d}}{\operatorname{d}t}
f(t)
=
\kappa
t
f(t)$
This condition can be interpreted as:
``derivative in time domain''
is equivalent to the
``'derivative in frequency domain''.
Therefore,
\begin{align*}
\frac{\operatorname{d}^2}{\operatorname{d}t^2}
f(t)
=
\kappa
\left[
f(t)
+ t
\cdot
\frac{\operatorname{d}}{\operatorname{d}t}
f(t)
\right]
=
\kappa f(t)
+
\kappa^2
t^2 f(t)
.
\end{align*}
Simple manipulations yield:
\begin{align*}
\frac{\operatorname{d}^2}{\operatorname{d}t^2}
f(t)
-
\kappa
(1+\kappa t^2)
(\kappa t)^2
\cdot
f(t)
=
0
.
\end{align*}
The only solutions on 
$\mathcal{E}$
correspond
to
$\kappa=\pm1$, i.e., 
$\frac{\operatorname{d}^2}{\operatorname{d}t^2}
f(t)
+
(1-t^2)f(t)=0$
or 
$\frac{\operatorname{d}^2}{\operatorname{d}t^2}
f(t)
-
(1+t^2)f(t)=0$.

\begin{proposition}
Any real signal
$f(t) \stackrel{\mathcal{F}}{\longleftrightarrow}  F(\omega)$
such that 
$f(t),
\frac{\operatorname{d}}{\operatorname{d}t}f(t),
F(\omega),
\frac{\operatorname{d}}{\operatorname{d}\omega}F(\omega)
\in L^2(\mathbb{R})$
have finite resolutions.
\end{proposition}
\proof
Applying the Parseval-Plancherel Theorem~\cite{ABRASTE68},
it follows that
\begin{align*}
\int_{-\infty}^\infty
t^2
f^2(t)
\operatorname{d}t
&
=
\int_{-\infty}^\infty
[j t f(t)]
\cdot
[j t f(t)]^\ast
\operatorname{d}t
\\
&
=
\frac{1}{2\pi}
\int_{-\infty}^\infty
\left|
\frac{\operatorname{d}}{\operatorname{d}\omega}
F(\omega)
\right|^2
\end{align*}
and
\begin{align*}
\int_{-\infty}^\infty
\omega^2
\left|
F^2(\omega)
\right|^2
\operatorname{d}\omega
&
=
\int_{-\infty}^\infty
[j \omega F(\omega)]
\cdot
[j \omega F(\omega)]^\ast
\operatorname{d}\omega
\\
&
=
2\pi
\int_{-\infty}^\infty
\left[
\frac{\operatorname{d}}{\operatorname{d}t}
f(t)
\right]^2
.
\end{align*}
Therefore,
\begin{align*}
\Delta_t^2
=
\frac{
\int_{-\infty}^\infty
\left|
\frac{\operatorname{d}}{\operatorname{d}\omega}
F(\omega)
\right|^2
\operatorname{d}\omega
}
{
\int_{-\infty}^\infty
\left|
F(\omega)
\right|^2
\operatorname{d}\omega
}
<
\infty
,
\\
\Delta_\omega^2
=
\frac{
\int_{-\infty}^\infty
\left|
\frac{\operatorname{d}}{\operatorname{d}t}
f(t)
\right|^2
\operatorname{d}t
}
{
\int_{-\infty}^\infty
\left|
f(t)
\right|^2
\operatorname{d}\omega
}
<
\infty
.
\end{align*}
Thus,
the above quantities
are
given by the square root of the ratio between the energy of the signal derivative and the energy of signal itself.
Thus,
the resolution for the Fourier invariant signal 
$\operatorname{sech}(\cdot)$
given by~\eqref{equation-4}
is
\begin{align*}
\Delta_t
=
\Delta_\omega
=
\sqrt{\frac{\pi}{6}}
\approx
0.7235987766\ldots
\end{align*}
since that
\begin{gather*}
\int_{-\infty}^\infty
\operatorname{sech}(t)
\operatorname{d}t
=2
,
\\
\int_{-\infty}^\infty
\operatorname{tanh}^2(t)
\operatorname{sech}^2(t)
\operatorname{d}t
=\frac{2}{3}
,
\\
\int_{-\infty}^\infty
\left(
\frac{2 t}{\pi}
\right)^2
\operatorname{sech}^2(t)
\operatorname{d}t
=\frac{2}{3}
,
\end{gather*}
where
$\operatorname{tanh}(\cdot)$
is the hyperbolic tangent function.

\begin{proposition}
[{\cite{GAB46}}]
Time-frequency uncertainty of Fourier Eigenfunctions
$\psi_n^\ast(t) = H_n(t) e^{-t^2/2} e^{j\omega_0 t + \phi_0}$,
where $\omega_0$ and $\phi_0$ are constants,
attain quantized values of the Gabor-Heisenberg lower bound, i.e.
\begin{align*}
\Delta t \cdot \Delta f = \frac{1}{2} \cdot (2n+1)
,
\\
\Delta_t \cdot \Delta_\omega = \frac{1}{2} \cdot (2n+1)
.
\end{align*}
\end{proposition}
That is why Gabor functions are 
relevant in some problems (e.g.~\cite{OKA98}.)

\section{The Concept of Isoresolution Wavelet}

The concept of isoresolution analysis is introduced in this section. According to the Gabor principle, if one increases resolution in one domain, the resolution must decrease in the other domain so as to guarantee the lower bound given by~\eqref{equation-20}. When analysing signals in joint time-frequency plane, frequently, there is no grounds to assure a better resolution in a domain than in the other domain. As an interesting property, any Fourier eigenfunction achieves isoresolution as it can be seen by the following proposition.

\begin{proposition}
Fourier-invariant signals perform an isoresolution, that is,
$\Delta t = \Delta_\omega$
\end{proposition}
\proof
Supposing that
$f\in\mathcal{E}$, 
then
$F(\omega) = \lambda f(\omega)$.
Therefore:
\begin{align*}
\frac{
\int_{-\infty}^\infty
F(\omega)
\omega^2
F(\omega)^\ast
\operatorname{d}\omega
}
{
\int_{-\infty}^\infty
\left|
F(\omega)
\right|^2
\operatorname{d}\omega
}
=
\frac{
\int_{-\infty}^\infty
\omega^2
|\lambda|^2
f^2(\omega)
\operatorname{d}\omega
}
{
\int_{-\infty}^\infty
\left|
\lambda
\right|^2
f^2(\omega)
\operatorname{d}\omega
}
\end{align*}
and the proof follows.
\endproof

\begin{table*}[t]
\centering

\caption{Resolution of a few standard continuous wavelets}
\label{table-1}

\begin{tabular}{lccc}
\toprule
Wavelet name
&
Time resolution $\Delta_t$
&
Frequency resolution $\Delta_\omega$
&
Isoresolution factor
$\sqrt{\Delta_t/\Delta_\omega}$
\\
\midrule
\texttt{Gaus1} & 1.500000 & 1.500000 & 1.000000
\\
\texttt{mexh}  & 1.166667 & 2.500000 & 0.683130
\\
\texttt{morl}  & 0.500002 & 25.499997 & 0.140028
\\
\texttt{fbsp 2-1-0.5} & $\infty$ & 14.475133 & -
\\
\texttt{shan 1-0.5} & $\infty$ & 13.159733 & -
\\
\texttt{haar} & 0.333333 & $\infty$
\\
\bottomrule
\end{tabular}
\end{table*}

This is an interesting property for signalling on the joint time-frequency plane. 

It is suggested here the changing of the time-frequency resolution by a proper scaling that allows for identical resolution in both domains.

\begin{proposition}
\label{proposition-9}
If $\psi(t)$
has effective duration
$\Delta t$
and
effective bandwidth 
$\Delta \omega$,
then the scaled version
$\psi\left( \sqrt{\Delta_t / \Delta_\omega}\,\, t \right)$
achieves isoresolution.
\end{proposition}
\proof
Scaled versions $\psi(at)$, $a\neq0$,
have resolutions
$\Delta_t/|a|$
and
$|a|\cdot \Delta_\omega$,
so $|a|$ can be appropriately chosen.
\endproof
The quantity
$\sqrt{\Delta_t / \Delta_\omega}$
is refered to as the \emph{isoresolution factor}.
The essential idea of isoresolution can be placed in the wavelet structure. Normally, the basic wavelet of a family 
$\frac{1}{\sqrt{|a|}}\psi\left(\frac{t-b}{a}\right)$
holds the admissibility condition but often 
does not achieve isoresolution.
We propose here to redefine 
the basic wavelet of a family so as to achieve isoresolution.
For instance,
the standard Mexican hat wavelet~$\psi_\text{Mhat}(t)$
satisfies:
\begin{align*}
2(t^2-1)
\cdot
\frac{e^{-t^2/2}}{\sqrt[4]{\pi}\sqrt{3}}
\stackrel{\mathcal{F}}{\longleftrightarrow}
-2
\sqrt{\frac{2}{3}}
\sqrt[4]{\pi}
\omega^2
e^{-\omega^2/2}
.
\end{align*}

The isoresolution Mexican hat wavelet can be found applying 
Proposition~\ref{proposition-9}:
\begin{align*}
\sqrt{\frac{7}{15}}
\cdot
\psi_\text{Mhat}
\left(
\sqrt{\frac{7}{15}}
\,\,
t
\right)
.
\end{align*}
For any isoresolution wavelet, 
the scaling by $a>1$ or $a<1$
corresponds to unbalance resolution in a different way.
Table~\ref{table-1}
displays both time and frequency resolution for a few known continuous wavelets: 
Gaussian derivatives, 
Mexican hat, 
Morlet, 
frequency B-Spline, 
Shannon, and 
Haar~\cite{MIS01}.
The wavelet \texttt{Gaus1} is an invariant wavelet therefore it achieves isoresolution, in accordance to proposition 8. It is valuable to mention that compact support wavelets (in time or frequency) cannot attain isoresolution, since no signal can simultaneously be 
time and frequency limited~\cite{WOZJAC67}.

\section{Perspectives and Closing Remarks}

Eigenfunctions of the Fourier operator were investigated and the Gabor principle was revisited defining the concept of isoresolution, i.e, a signal with the same time and frequency resolution.
The  functions $\left\{ \psi_n(t) \right\}$
(see~\eqref{equation-15}) turn up as a very appealing choice for designing representations such as wavelets. It is time to try finding new wavelets starting with~\eqref{equation-11}. Since they are solutions of a wave equation (2nd order differential equation), our approach 
(Mathieu~\cite{LIRA03a}, Legendre~\cite{LIRA03b}, 
Chebyshev~\cite{CIN03}) can be useful to construct new wavelets: The Quantum Wavelets, or Gabor-Schr\"odinger wavelets. The construction of new wavelets based on these complete, orthogonal, domain shape-invariant system is currently being investigated. The idea is to adapt the concept of isoresolution in orthogonal multiresolution analysis~\cite{DEO03, HESSWIK96}. 

\section*{Acknowledgments}

This work was partially supported by CNPq and CAPES.
The authors also thank Mr.~M.~M\"uller for some motivating discussions.

{\small
\bibliographystyle{IEEEtran}
\bibliography{ref-clean}

\begin{thebibliography}{10}
\providecommand{\url}[1]{#1}
\csname url@samestyle\endcsname
\providecommand{\newblock}{\relax}
\providecommand{\bibinfo}[2]{#2}
\providecommand{\BIBentrySTDinterwordspacing}{\spaceskip=0pt\relax}
\providecommand{\BIBentryALTinterwordstretchfactor}{4}
\providecommand{\BIBentryALTinterwordspacing}{\spaceskip=\fontdimen2\font plus
\BIBentryALTinterwordstretchfactor\fontdimen3\font minus
  \fontdimen4\font\relax}
\providecommand{\BIBforeignlanguage}[2]{{%
\expandafter\ifx\csname l@#1\endcsname\relax
\typeout{** WARNING: IEEEtran.bst: No hyphenation pattern has been}%
\typeout{** loaded for the language `#1'. Using the pattern for}%
\typeout{** the default language instead.}%
\else
\language=\csname l@#1\endcsname
\fi
#2}}
\providecommand{\BIBdecl}{\relax}
\BIBdecl

\bibitem{HERS64}
I.~N. Herstein, \emph{Topics in Algebra}.\hskip 1em plus 0.5em minus
  0.4em\relax Mass.: Blaisdell Pub. Co., 1964.

\bibitem{SOKORED66}
I.~S. Sokolnikoff and R.~M. Redheffer, \emph{Mathematics of Physics and Modern
  Engineering}, 2nd~ed.\hskip 1em plus 0.5em minus 0.4em\relax Tosho:
  McGraw-Hill Kogakusha, 1966.

\bibitem{PEIDIN02}
S.~C. Pei and J.~J. Ding, ``Eigenfunctions of linear canonical transform,''
  \emph{IEEE Trans. on Signal Processing}, vol.~50, no.~1, pp. 11--26, Jan.
  2002.

\bibitem{COH95}
A.~Cohen, \emph{Time-Frequency Analysis}.\hskip 1em plus 0.5em minus
  0.4em\relax Prentice-Hall Signal Processing Series, 1995.

\bibitem{QIACHE99}
S.~Qian and D.~Chen, ``Joint time-frequency analysis,'' \emph{IEEE Signal Proc.
  Mag.}, vol.~16, no.~2, pp. 52--67, Mar. 1999.

\bibitem{ABRASTE68}
M.~Abramowitz and I.~Stegun, Eds., \emph{Handbook of Mathematical
  Functions}.\hskip 1em plus 0.5em minus 0.4em\relax New York: Dover, 1968.

\bibitem{BEI67}
A.~Beiser, \emph{Concepts of Modern Physics}.\hskip 1em plus 0.5em minus
  0.4em\relax McGraw-Hill Series in Fundamental Physics, 1967.

\bibitem{OLIBAR00}
P.~M. Oliveira and V.~Barroso, ``Uncertainty in the time-frequency plane,'' in
  \emph{Proc. of the Tenth IEEE Workshop on Statistical Signal and Array
  Processing}, Pocono Manor, PA , USA, 2000, pp. 607--611.

\bibitem{GAB46}
D.~Gabor, ``Theory of communications,'' \emph{J. IEE}, vol.~93, pp. 429--457,
  1946.

\bibitem{GAB53}
------, ``Communication theory and physics,'' \emph{IEEE Trans. Info. Theory},
  vol.~1, pp. 48--59, Feb. 1953.

\bibitem{OKA98}
K.~Okajima, ``The {G}abor function extract the maximum information from input
  local signals,'' \emph{Neural Networks}, vol.~11, pp. 435--439, 1998.

\bibitem{MIS01}
M.~Misiti, M.~Misiti, G.~Oppenheim, and J.~M. Poggi, \emph{Wavelet
  Toolbox}.\hskip 1em plus 0.5em minus 0.4em\relax The Math Works, 2001.

\bibitem{WOZJAC67}
J.~M. Wozencraft and I.~M. Jacobs, \emph{Principles of Communication
  Engineering}.\hskip 1em plus 0.5em minus 0.4em\relax New York: Wiley, 1967.

\bibitem{LIRA03a}
\BIBentryALTinterwordspacing
M.~M.~S. Lira, H.~M. Oliveira, and R.~J. Cintra, ``Elliptic-cylindrical
  wavelets: the {M}athieu wavelets,'' \emph{IEEE Signal Processing Letters},
  vol.~11, no.~1, pp. 52--55, Jan. 2004. [Online]. Available:
  \url{http://ieeexplore.ieee.org/xpl/freeabs_all.jsp?arnumber=1255923}
\BIBentrySTDinterwordspacing

\bibitem{LIRA03b}
\BIBentryALTinterwordspacing
M.~M.~S. Lira, H.~M. de~Oliveira, M.~A. {Carvalho Jr}, and R.~M. {Campello de
  Souza}. Compactly supported wavelets derived from {L}egendre polynomials:
  Spherical harmonic wavelets. ArXiv:1502.00950. [Online]. Available:
  \url{http://arxiv.org/abs/1502.00950}
\BIBentrySTDinterwordspacing

\bibitem{CIN03}
\BIBentryALTinterwordspacing
R.~J. Cintra, L.~R. Soares, and H.~M. de~Oliveira. Filter banks and wavelets
  based on {C}hebyshev polynomials. ArXiv:1411.2389. [Online]. Available:
  \url{http://arxiv.org/abs/1411.2389}
\BIBentrySTDinterwordspacing

\bibitem{DEO03}
H.~M. de~Oliveira, L.~R. Soares, and T.~H. Falk, ``A family of wavelets and a
  new orthogonal multiresolution analysis based on the {N}yquist criterion,''
  \emph{Rev. da Soc. Bras. Telecomm.}, Jun. 2003.

\bibitem{HESSWIK96}
N.~Hess-Nielsen and M.~V. Winckerhauser, ``Wavelets and time-frequency
  analysis,'' \emph{Proc. of the IEEE}, vol.~84, no.~4, pp. 523--540, Apr.
  1996.

\end{thebibliography}
}

\end{document}